# DISJUNCTIVE IDENTITIES OF FINITE GROUPS AND IDENTITIES OF REGULAR REPRESENTATIONS

Samuel M. Vovsi

Institute for Advanced Study

## 1. Introduction

The problem of describing identities of concrete algebraic structures is traditionally one of the most attractive in variety theory. Very often it is a difficult problem, and many deep works are entirely devoted to finding bases of identities of classical algebraic objects.

This paper was initiated by attempts to find bases of identities for certain representations of finite groups. Let $F$ be the free group of countable rank on an alphabet $x_1, x_2, \ldots$, and let $KF$ be its group algebra over a field $K$. Recall that an element $u = u(x_1, \ldots, x_n) \in KF$ is called an *identity* of a representation $\rho : G \to \operatorname{Aut}_K V$ (usually denoted by $\rho = (V, G)$) if for arbitrary $g_1, \ldots, g_n \in G$ the element $u(g_1, \ldots, g_n) \in KG$ annihilates $V$. It is well known that if $\rho$ is a representation of a *finite* group, then its identities are *finitely based*, that is, are in a natural sense generated by some finite set of them [V, VN]. The problem is how to find such a finite basis *explicitly*. Of course, there is no hope to develop a universal method allowing to write down a basis of identities for an arbitrary representation, so we will be concerned with a much more concrete (but still very general)

**Problem A.** *Let $G$ be a finite group, $\operatorname{char} K \nmid |G|$. Find a basis of identities for the regular representation $\operatorname{Reg}_K G = (KG, G)$ of $G$ over $K$.*

Why the regular representations are of particular interest? The reason is two-fold. First, it is easy to understand that the identities of $\operatorname{Reg}_K G$ are the identities of *all* representations of $G$ over $K$. Second, they can be regarded as special "weak" identities of the group algebra $KG$, for they are exactly the elements $u(x_1, \ldots, x_n) \in KF$ such that $u(g_1, \ldots, g_n) = 0$ in $KG$ for all $g_i \in G$.

An interesting approach to Problem A was suggested by Plotkin and Kushkuley [PK]. Consider a universally quantified formula

$$\delta = \forall x_1 \ldots x_n[(f_1 = 1) \vee (f_2 = 1) \vee \cdots \vee (f_n = 1)] \tag{1}$$

where $x_1, \ldots, x_n$ are all variables involved in the $f_i$. Probably, the first systematic study of such formulas was done by Baker [B], who called them universal disjunctions of equations

Typeset by $\mathcal{A}\mathcal{M}\mathcal{S}$-TEX

(UDE). We say that $\delta$ is a *disjunctive identity* (d-identity, for short) of a group $G$ if it is identically true on $G$. It is obvious that if $\delta$ is a d-identity of $G$ then

$$\delta^* = (f_1^{y_1} - 1)\ldots(f_n^{y_n} - 1) \in KF$$

is an identity of every representation of $G$. A principal observation made in [PK] is that if one knows the d-identities of a finite group $G$, then one can easily obtain a basis of identities of the regular representation $\text{Reg}_K G$. To make this more precise, we recall a few simple definitions and facts.

Let $\delta$ be an arbitrary d-identity, that is, a universal sentence of the form (1). Usually we omit the quantor prefix and write this formula simply as $\delta = (f_1 = 1) \vee (f_2 = 1) \vee \cdots \vee (f_n = 1)$. A class of groups $\mathfrak{D}$ is called a *d-variety* if it is definable by a set of d-identities.[1] This is equivalent to saying that $\mathfrak{D}$ is definable by positive universal sentences, because every such a sentence is obviously equivalent to a conjunction of d-identities. Applying a well known fact of model theory [G; p. 275], we see that a class of groups is a d-variety if and only if it is closed under taking subgroups, homomorphic images, and ultraproducts.

If $\mathfrak{X}$ is an arbitrary class of groups, then by $\text{dvar}(\mathfrak{X})$ we denote the d-variety generated by $\mathfrak{X}$. In other words, $\text{dvar}(\mathfrak{X})$ is the class of all groups satisfying all d-identities that are valid in all groups from $\mathfrak{X}$. Let $G$ be a finite group and let $\text{HS}(G)$ be the class of all sections of $G$. Being a finite class of finite algebras, $\text{HS}(G)$ is axiomatizable by a single sentence. Therefore it is closed under ultraproducts and so is a d-variety. Thus we have established

**Fact 1.** *The d-variety generated by a finite group is the class of all sections of this group.*[2]

In particular, this implies that if $G$ and $H$ are finite groups then

$$\text{dvar}(G) = \text{dvar}(H) \iff G \simeq H.$$

In other words, *every finite group is uniquely determined by its d-identities.*

A set of d-identities $\{\delta_i\}$ of a class of groups $\mathfrak{X}$ is called a *basis of d-identities of $\mathfrak{X}$* if every group satisfying the $\delta_i$ belongs to $\text{dvar}(\mathfrak{X})$ (that is, satisfies *all* d-identities of $\mathfrak{X}$). It is easy to see that every finite group has a finite basis of d-identities. Indeed, let $G$ be group of order $n$. Consider the d-identity

$$\omega_n = \bigvee_{0 \leq i < j \leq n} (x_i = x_j)$$

and for every group $H_i$ of order at most $n$ and not belonging to $\text{dvar}(G)$ choose a d-identity $\delta_i$ of $G$ which is not valid in $H_i$. Then it is clear that a finite set of formulas $\{\omega_n, \delta_1, \delta_2, \ldots\}$ is a basis of d-identities of $G$.[3]

---

[1] Initially, Plotkin and Kushkuley [PK] used the words "pseudoidentity" and "pseudovariety". We would also prefer to use this (pseudo)terminology but, unfortunately, it has been already occupied by Eilenberg and Schützenberger [ES].

[2] A direct proof of this fact can be found, for example, in [PV; Theorem 13.2.1].

[3] All of the above is based on general results of model theory and is valid not only for groups, but for arbitrary algebraic systems

We now discuss the connection between d-identities of groups and identities of their regular representations. We will consider representations of groups over an arbitrary but fixed field $K$. As defined earlier, identities of representations are elements of the group algebra $KF$. A class of group representations $\mathcal{V}$ is called a *variety* if it consists of all representations satisfying a given set of identities. A variety is *finitely based* if it is definable by a finite set of identities. Such a set is called a (finite) *basis* of identities of $\mathcal{V}$. If $\mathcal{V}$ is the variety generated by some class of representations $\mathcal{X}$ (notation: $\mathcal{V} = \text{var}(\mathcal{X})$), then the same set of identities is also called a basis of identities of $\mathcal{X}$.

Let $\mathfrak{V}$ be a variety of groups. Take the class of all representations $\rho : G \to \text{Aut}_K V$ such that $G/\text{Ker}\,\rho \in \mathfrak{V}$. It is a variety of group representations: if $\mathfrak{V}$ is definable by group identities $f_i$, then $\omega\mathfrak{V}$ is defined by the identities $f_i - 1$. A variety of group representations $\mathcal{V}$ is *locally finite* if it satisfies the following condition: if $\rho = (V, G)$ is an arbitrary representation in $\mathcal{V}$ such that $G$ is a finitely generated group and $V$ is a finitely generated $G$-module, then $V$ is finite-dimensional over $K$ and the group $G/\text{Ker}\,\rho$ is finite. This is equivalent to saying that there exists a locally finite variety of groups $\mathfrak{V}$ such that $\mathcal{V} \subseteq \omega\mathfrak{V}$ (see, for example, [V2; §2.1]) . If in addition this $\mathfrak{V}$ can be chosen in such a way that $\exp \mathfrak{V}$ is not divisible by the characteristic of the ground field $K$, the variety $\mathcal{V}$ is called *ordinary*, or *nonmodular*.

**Example.** Let $\rho = (V, G)$ be an ordinary representation of a finite group $G$ (that is, $\text{char}\,K \nmid |G|$). Then it is easy to see that the variety generated by $\rho$ is locally finite and ordinary. It is known that this variety is finitely based [V1]. Moreover, the latter is true for arbitrary (not necessary ordinary) representations of finite groups [VN].

Let $G$ be a group and let $\{\delta_i\}$ be a set of d-identities of $G$. This set is called a *weak basis* of d-identities of $G$ if every irreducible linear (over $K$) group satisfying all the $\delta_i$ belongs to $\text{dvar}(G)$ (i.e. is a section of $G$). Of course, every basis of d-identities is a weak basis, but the converse is not, in general, true. The following fact plays a key role.

**Fact 2** [PK]. *Let $G$ be a finite group and let $U$ be a set of identities of the regular representation $\text{Reg}_K G = (KG, G)$, satisfying the following conditions:*
   (i) *the variety of representations defined by $U$ is locally finite and ordinary;*
   (ii) *there exists a weak basis $\{\delta_i\}$ of d-identities of $G$ such that $U \supseteq \{\delta_i^*\}$.*
*Then $U$ is a basis of identities of $\text{Reg}_K G$.*

P r o o f  of this statement is not difficult, and we will provide it here. Consider a d-identity $\delta = (f_1 - 1) \vee \cdots \vee (f_n - 1)$, and let $\delta^* = (f_1^{y_1} - 1) \ldots (f_n^{y_n} - 1)$. First let us show that if $\delta^*$ is an identity of some faithful irreducible representation $\rho = (V, H)$, then $\delta$ must be an identity of the acting group $H$. Assume the contrary; then there is a homomorphism $\phi : F \to H$ such that $f_i^\phi \neq 1$ for every $i$. Denote $f_i^\phi = g_i$ and take $0 \neq a_1 \in V$. Since $V$ is an irreducible $KH$-module, we have $a_1 KH = V$. Since $H$ acts on $V$ faithfully and $g_1 \neq 1$, there exists $h_1 \in H$ such that $a_1 \cdot h_1^{-1} g_1 \neq a_1 \cdot h_1^{-1}$, that is, $a_1 \cdot (g_1^{h_1} - 1) \neq 0$. Next, denote $a_1 \cdot (g_1^{h_1} - 1) = a_2$ and find $h_2 \in H_2$ such that $a_3 = a_2 \cdot (g_2^{h_2} - 1) \neq 0$. Repeating this argument, we find elements $h_3, \ldots, h_n \in H$ such that
$$a \cdot (g_1^{h_1} - 1)(g_2^{h_2} - 1) \ldots (g_n^{h_n} - 1) \neq 0.$$

But this is impossible because $\delta^*$ is an identity of $\rho$.

Now let $G$ and $U$ be as in the statement. Let $\mathcal{V}$ be the variety generated by $\operatorname{Reg} G$ and let $\mathcal{W}$ be the variety defined by $U$. Then $\mathcal{V} \subseteq \mathcal{W}$; we have to show that $\mathcal{V} = \mathcal{W}$. By (ii), $\mathcal{W}$ is locally finite and ordinary, and therefore it is generated by faithful irreducible representations of finite groups. Take in $\mathcal{W}$ such a representation $\rho = (V, H)$; it remains to show that $\rho \in \mathcal{V}$.

Since the $\delta_i^*$ belong to $U$, they all are identities of $\rho$. By the above, the $\delta_i$ are d-identities of $H$. Since $H$ is an irreducible linear group and $\{\delta_i\}$ is a weak basis of d-identities of $G$, it follows that $H \in \operatorname{dvar}(G)$, that is, $H$ is a section of $G$. But then $\operatorname{Reg} H$ satisfies all identities that are valid in $\operatorname{Reg} G$. Hence $\operatorname{Reg} H$ belongs to $\mathcal{V}$, and the same is true for *any* representation of $H$. □

Since it is very easy to find identities of the representation $\operatorname{Reg}_K G$ which guarantee that the corresponding variety is locally finite and ordinary, the established fact reduces solution of Problem A to finding a weak basis of d-identities of a given finite group $G$. These ideas were developed in [PK] in several directions; in particular, they were applied to the study of identities of the regular representations of the symmetric groups $S_n$, $n \leq 5$.

The present paper is focused on finding *actual* (not weak) bases of d-identities.

**Problem B.** *For a given finite group $G$, find a basis of its d-identities.*

This problem is more difficult than its "weak sibling" and is of much greater independent interest. Indeed, as mentioned earlier, every finite algebra is uniquely determined by its disjunctive identities, and therefore a basis of d-identities of a finite group uniquely identifies it. What is especially important, Problem B is in a certain sense "finite", and if the order of the group in question is small enough, it can be solved by a rather mechanical procedure.

The aim of this paper is to collect a sufficient experimental material. It consists of six sections. In §2 we compute bases of d-identities for several small groups (like $D_8, Q_8, S_4$), but in §3 — bases of d-identities for $A_5$ and $S_5$. In particular, in §3 we sharpen some results from [PK] and eliminate a few gaps from that work. §4 is devoted to finding a basis of d-identities for $A_6$; technically it is the most difficult part of the paper. In §5, contrary to the previous sections, we find bases of d-identities for an *infinite series* of groups, namely, for all dihedral groups $D_{2n}$. The results of this section are joint with Gregory Cherlin. Finally, in §6 we return to our initial Problem A and explicitly write down bases of identities for the (nonmodular) regular representations of certain groups discussed in Sections 2–5.

The results of this work demonstrate some general principles of finding bases for d-identities of finite groups. It is not excluded that a certain computer-based system (like GAP or CAYLEY) can be involved in this process, and we hope to investigate such a possibility in the future.

It should be noted that we essentially use some ideas and results from [PK]. Unfortunately, this interesting work has never been published in full form and exists only as a preprint (part of it was expounded in [PV; §13]). We are grateful to B. I. Plotkin for making this preprint available.

## 2. Several simple examples

The aim of this section is to illustrate, on several simple examples, some basic principles of finding bases of d-identities of finite groups. We will consider the groups $D_8$, $Q_8$, $A_4$, and $S_4$.

THE GROUP $D_8$. The dihedral group of order eight has the following presentation:

$$D_8 = \langle a, b \,|\, a^4 = b^2 = 1, a^b = a^{-1} \rangle.$$

Clearly it satisfies the d-identities

$$\omega_8 = \bigvee_{0 \leq i < j \leq 8} (x_i = x_j) \qquad (2.1)$$

and

$$x^4 = 1. \qquad (2.2)$$

Since $D_8$ contains only two elements of order $> 2$, it also satisfies the d-identity

$$(x_1^2 = 1) \vee (x_2^2 = 1) \vee (x_3^2 = 1) \vee (x_1 = x_2) \vee (x_1 = x_3) \vee (x_2 = x_3). \qquad (2.3)$$

Let us prove that the formulas (2.1)–(2.3) form a *weak* basis of d-identities of $D_8$. We have to show that every irreducible linear group $G$ satisfying (2.1)–(2.3) belongs to $\mathrm{dvar}(D_8)$, i.e. is a section of $D_8$. Since $G$ satisfies $\omega_8$, it must be of order at most 8. The groups of orders 3,5,6,7 are excluded by (2.2), but the groups of orders 1,2,4 are all sections of $D_8$. There are five groups of order 8:

$$Z_8,\ Z_4 \times Z_2,\ Z_2 \times Z_2 \times Z_2,\ Q_8,\ \text{and}\ D_8,$$

where, as usual, $Z_n$ denotes the cyclic group of order $n$ and $Q_8$ the quaternion group. $Z_8$ does not satisfy (2.2), each of the groups $Z_4 \times Z_2$ and $Q_8$ has three distinct elements of order 4 and so does not satisfy (2.3), and $Z_2^3$ is not irreducible linear (an abelian irreducible linear group must be cyclic). This completes the proof.

However, the formulas (2.1)–(2.2) do not form a basis of d-identities of $D_8$ because $Z_2^3$ obviously satisfies all of them but is not a section of $D_8$. To eliminate $Z_2^3$, we have to introduce another d-identity. For brevity, for any variables $x$ and $y$ we will write $x \in \langle y \rangle$ to denote

$$(x = 1) \vee (x = y) \vee (x = y^2) \vee (x = y^3).$$

Consider the following formula in the variables $x_1, x_2, x_3$:

$$\Big[\bigvee_{i \neq j} x_i \in \langle x_j \rangle\Big] \vee \Big[\bigvee_{\sigma \in S_3} x_{\sigma(1)} \in \langle x_{\sigma(2)} x_{\sigma(3)} \rangle\Big] \vee \Big[\bigvee_{\sigma \in S_3} x_{\sigma(1)} x_{\sigma(2)} \in \langle x_{\sigma(2)} x_{\sigma(3)} \rangle\Big] \qquad (2.4)$$

and show that it is a d-identity of $D_8$. The group $D_8$ is a semidirect product

$$D_8 = A \rtimes B, \quad \text{where} \quad A \simeq Z_4, \ B \simeq Z_2.$$

Note that if $x, y$ are elements of $D_8$ not in $A$, then $xy \in A$. Now take arbitrary $x_1, x_2, x_3 \in D_8$. If two of the $x_i$ lie in $A$ then the first clause of (2.4) holds, if one of them lies in $A$ and the other two outside then the second clause holds, and if all three $x_i$ lie outside then the third clause holds.

On the other hand, (2.4) is not valid in $Z_2^3$: it is enough to take a basis of this group for $x_1, x_2, x_3$. Thus we have proved

**Proposition 1.** *Formulas (2.1)–(2.4) form a basis of d-identities of $D_8$.* $\square$

Although the above arguments are rather elementary, they give a clear idea about the procedure of finding a basis of d-identities of a finite group $G$. If $|G| = n$, then we first take $\omega_n$, and then for every group $H$ of order $\leq n$ find a d-identity of $G$ that is not valid in $H$. It is clear that the complexity of this procedure increases drammatically with the order of the group in question. In general, it is noticeably easier to find a *weak* basis of d-identities because in this case we do not have to think about groups not admitting faithful irreducible representations. In particular, we do not have to think about groups of the form $Z_p^n$ which usually create most of the troubles.

THE GROUP $Q_8$. Recall that

$$Q_8 = \{\pm 1, \pm i, \pm j, \pm k\},$$

where

$$i^2 = j^2 = k^2 = -1,$$

$$ij = -ji = k, \quad jk = -kj = i, \quad ki = -ik = j,$$

and the other multiplications are natural. Obviously $Q_8$ satisfies (2.1) and (2.2). Every group satisfying these two d-identities has order 1, 2, 4, or 8. The groups of order 1, 2 and 4 are sections of $Q_8$, but $Z_8$ does not satisfy (2.2). It remains to eliminate three groups of order 8: $Z_4 \times Z_2$, $Z_2^3$ and $D_8$. Note that

$$(x_1 \in \langle x_2 \rangle) \vee (x_2 \in \langle x_1 \rangle) \vee (x_1^2 x_2^2 = 1), \tag{2.5}$$

where $x \in \langle y \rangle$ has the same meaning as before, is valid in $Q_8$. For if $x_1 \notin \langle x_2 \rangle$ and $x_2 \notin \langle x_1 \rangle$, then both $x_i$ are different from $\pm 1$ and therefore are of order 4. Hence $x_1^2 = x_2^2 = -1$ and $x_1^2 x_2^2 = 1$.

Next we observe that (2.5) eliminates both $Z_4 \times Z_2$ and $D_8$. Indeed in the first case we simply take a basis of $Z_4 \times Z_2$ for $x_1$ and $x_2$, but in the second case let $x_1$ and $x_2$ be the standard generators of $D_8$ (with $|x_1| = 4$ and $|x_2| = 2$). Then $x_1 \notin \langle x_2 \rangle$ and $x_2 \notin \langle x_1 \rangle$ and $x_1^2 x_2^2 \neq 1$.

To eliminate $Z_2^3$, consider the following formula in $x_1, x_2, x_3$:

$$(\bigvee_{i \neq j} x_i \in \langle x_j \rangle) \vee (x_1 x_2 = x_3) \vee (x_1 x_2 = x_3^{-1}). \tag{2.6}$$

Choosing a basis of $Z_2^3$ for $x_1, x_2, x_3$, we see that (2.6) is not valid in this group. Now let $x_1, x_2, x_3 \in Q_8$. If one of the $x_i$ is equal to $\pm 1$, then the first clause of (2.6) holds, otherwise all $x_i$ are elements of order 4. There are three cyclic subgroups of order 4 in $Q_8$, generated by $i$, $j$ and $k$, respectively. If two of the $x_i$ belong to one of these subgroups then again the first clause holds, but if all three of them lie in different subgroups then either the second or the third clause holds. Thus (2.6) is a d-identity of $Q_8$, and we obtain

**Proposition 2.** *Formulas (2.1), (2.2), (2.5) and (2.6) form a basis of d-identities of* $Q_8$. □

THE GROUP $A_4$. The alternating group of degree 4 contains three involutions (12)(34), (13)(24), (14)(23); all other nonunit elements are cycles of length 3. Therefore $A_4$ satisfies the d-identity

$$(x^2 = 1) \vee (x^3 = 1) \tag{2.7}$$

and, of course, $\omega_{12}$. The order of a group satisfying $\omega_{12}$ and (2.7) must be one of the numbers

$$1,\ 2,\ 3,\ 4,\ 6,\ 8,\ 9,\ 12.$$

All groups of orders 1,2,3 are sections of $A_4$. The cyclic group of order 4 does not satisfy (2.7), but the noncyclic one is a subgroup of $A_4$. There are two groups of order 6: $Z_6$ and $S_3$. The first one does not satisfy (2.7); as to the second one, consider the formula

$$(x^3 = 1) \vee (x_2^3 = 1) \vee ((x_1 x_2)^2 = 1). \tag{2.8}$$

Seting $x_1 = (12)$, $x_2 = (13)$, we see that it is not valid in $S_3$. But it is valid in $A_4$. Indeed, take arbitrary $x_1, x_2 \in A_4$. Observe that

$$H = \{1,\ (12)(34),\ (13)(24),\ (14)(23)\}$$

is a normal subgroup of exponent 2 in $A_4$. If one of the $x_i$ has order 3, we are done; otherwise both $x_1$ and $x_2$ belong to $H$, hence $x_1 x_2 \in H$ and $(x_1 x_2)^2 = 1$.

All groups of order 8 except $Z_2^3$ have elements of order 4 and therefore do not satisfy (2.7). To eliminate $Z_2^3$, we note that the unique Sylow 2-subgroup $H$ of $A_4$ is a 2-dimensional vector space over $GF(2)$, while $Z_2^3$ is 3-dimensional over $GF(2)$. Therefore the d-identity

$$\bigvee_{(k_1, k_2, k_3)} (x_1^{k_1} x_2^{k_2} x_3^{k_3} = 1),$$

where $(k_1, k_2, k_3)$ runs over all nonzero $\{0, 1\}$-vectors of length 3, is valid in $H$ but not in $Z_2^3$. Further, for every $x \in A_4$ we have $x^3 \in H$. It follows that

$$\bigvee_{(k_1, k_2, k_3)} (x_1^{3k_1} x_2^{3k_2} x_3^{3k_3} = 1) \tag{2.9}$$

is a d-identity of $A_4$ but not of $Z_2^3$.

To eliminate the groups of order 9, we first prove a simple lemma which will be repeatedly used in the sequel.

**Lemma 1.** *Let $\alpha$ and $\beta$ be two n-cycles in $S_{n+1}$. Then for some $d \leq n - 1$ the permutation $\alpha\beta^d$ is not an n-cycle.*

P r o o f. a) Suppose that $\mathrm{Supp}(\alpha) = \mathrm{Supp}(\beta)$. Without loss of generality, let

$$\alpha = (1, 2, 3, \ldots, n), \quad \beta = (1, b_2, b_3, \ldots, b_n)$$

where $2 \leq b_i \leq n$. For any $d$, the permutation $\alpha\beta^d$ fixes $n+1$, hence it is enough to show that for some natural $d \leq n-1$ it also fixes the number 1. Let $d = d_\beta(2,1)$ be the *distance from 2 to 1 in $\beta$*, that is, $d$ is the smallest number with $2(\beta^d) = 1$. Then $d \leq n-1$ and $1(\alpha\beta^d) = (1\alpha)\beta^d = 2\beta^d = 1$, as required.

b) Suppose that $\mathrm{Supp}(\alpha) \neq \mathrm{Supp}(\beta)$. We may assume that

$$\alpha = (1, 2, 3, \ldots, n), \quad \beta = (b_1, b_2, \ldots, b_{n-1}, n+1)$$

where $1 \leq b_i \leq n-1$. Note that $(n-1)\alpha\beta^d = n$ for any $d$. Set $d = d_\beta(1, n-1)$, then $n(\alpha\beta^d) = (n\alpha)\beta^d = 1\beta^d = n-1$ and so the cycle form of $\alpha\beta^d$ contains a transposition $(n-1, n)$. It follows that $\alpha\beta^d$ is not an $n$-cycle. $\square$

It is now easy to see that

$$(x_1^4 = 1) \vee (x_2^4 = 1) \vee ((x_1x_2)^4 = 1) \vee ((x_1x_2^2)^4 = 1) \tag{2.10}$$

is a d-identity of $S_4$ (and so of $A_4$). Indeed, if $x_1, x_2 \in S_4$ and $x_i^4 \neq 1$, then $x_1, x_2$ are 3-cycles. By Lemma 1, either $x_1x_2$ or $x_1x_2^2$ is not a 3-cycle, and (2.10) follows. On the other hand, (2.10) can not be valid in a group of order 9.

Finally, all groups of order 12 except $A_4$ contain elements of order 6 and do not satisfy (2.7). Puting all this together, we have

**Proposition 3.** *Formulas $\omega_{12}$ and (2.7)–(2.10) form a basis of d-identities of $A_4$.* $\square$

Note that (2.9) was only used to eliminate $Z_2^3$, but (2.10) — to eliminate $Z_3^2$ (the only other group of order 9 is cyclic and does not satisfy (2.7)). Since none of these two groups admits a faithful irreducible representation, it follows that the formulas $\omega_{12}$, (2.7) and (2.8) form a *weak* basis of d-identities of $A_4$.

THE GROUP $S_4$. First we introduce three new d-identities:

$$(x^3 = 1) \vee (x^4 = 1), \tag{2.11}$$

$$(x_1^6 = 1) \vee (x_2^6 = 1) \vee (x_1 = x_2) \vee (x_1x_2 = 1) \vee ((x_1x_2)^3 = 1), \tag{2.12}$$

$$(x_1^3 = 1) \vee (x_2^3 = 1) \vee (x_3^3 = 1) \vee$$
$$\vee ((x_1^{-1}x_2)^3 = 1) \vee ((x_1^{-1}x_3)^3 = 1) \vee ((x_2^{-1}x_3)^3 = 1) \vee \theta, \quad (2.13)$$

where $\theta$ denotes the formula (2.4).

**Proposition 4.** *Formulas $\omega_{24}$ and (2.10)–(2.13) form a basis of d-identities of $S_4$.*

P r o o f. a) First we show that these d-identities are valid in $S_4$. This was already verified for (2.10) and is obvious for $\omega_{24}$ and (2.11). Next, take $x_1, x_2 \in S_4$ such that $x_i^6 \neq 1$. Then the $x_i$ are 4-cycles and a straightforward verification shows that, unless $x_1 = x_2^{\pm 1}$, $x_1 x_2$ is a 3-cycle. Hence (2.12) is valid in $S_4$.

Further, take arbitrary $x_1, x_2, x_3 \in S_4$ and show that they satisfy (2.13). If one of the $x_i$ is a 3-cycle, we are done; otherwise all the $x_i$ are 2-elements. Suppose they all belong to the same Sylow 2-subgroup of $S_4$. This subgroup is isomorphic to $D_8$; since $\theta$ is identically true in $D_8$, the $x_i$ satisfy (2.13). Now suppose that two of the $x_i$, say $x_1$ and $x_2$, belong to distinct Sylow 2-subgroups of $S_4$. The Sylow 2-subgroup of $S_4$ are

$$H_1 = \langle H, (12) \rangle, \quad H_2 = \langle H, (13) \rangle, \quad H_3 = \langle H, (23) \rangle,$$

where $H$ is as above. Let for example $x_1^{-1} = h_1 \, (12)$, $x_2 = h_2 \, (13)$ with $h_i \in H$. Then $x_1^{-1}x_2 = h_3 \, (123)$ with $h_3 \in H$. If $(x_1^{-1}x_2)^3 \neq 1$, then $x_1^{-1}x_2$ lie in some Sylow 2-subgroup of $S_4$, say $x_1^{-1}x_2 = h_4 \, (23)$ with $h_4 \in H$. But then $(123)(23) = h_3^{-1}h_4 \in H$ which is impossible because $(123)(23)$ is an odd permutation. Thus $(x_1^{-1}x_2)^3 = 1$ and (2.13) holds.

b) Now we show that every group satisfying $\omega_{24}$ and (2.10)–(2.13) is a section of $S_4$. The order of a group satisfying $\omega_{24}$ and (2.11) must be among the numbers

$$1, 2, 3, 4, 6, 8, 9, 12, 16, 18, 24.$$

All groups of order $\leq 4$ are subgroups of $S_4$. The cyclic group of order 6 does not satisfy (2.11) but the noncyclic one is a subgroup of $S_4$. Among the five groups of order 8, $Z_8$ does not satisfy (2.11); $Z_4 \times Z_2$ does not satisfy (2.12) (set $x_1 = (a,1)$, $x_2 = (a,b)$ where $a, b$ is the natural basis of $Z_4 \times Z_2$); $Z_2^3$ is eliminated by (2.13) (take a basis of this group for the $x_i$); $D_8$ is a subgroup of $S_4$; and $Q_8$ is eliminated by (2.12) (since there are two elements $i$ and $j$ of order 4 in $Q_8$ whose product $ij = k$ is also of order 4).

Every group of order 9 does not satisfy (2.10). The same is true for groups of order 18. Among the groups of order 12, $A_4$ is a subgroup of $S_4$ but the others contain elements of order 6 and do not satisfy (2.11). It is easy to see that every group of order 16 has an abelian section of order 8, and therefore is eliminated by either (2.11), or (2.12), or (2.13). It remains to consider groups of order 24.

**Lemma 2.** *Let $G$ be a group of order 24 without elements of order 6. Then $G \simeq S_4$.*

P r o o f. Suppose there exists $A \triangleleft G$ such that $|A| = 3$. Take in $A$ an element $a$ of order 3, and let $S$ be a Sylow 2-subgroup of $G$. Then $|S/C_S(A)|$ is 1 or 2. In both cases $S$ has an element $b$ of order 2 centralizing $A$. But then $ab$ is an element of order 6, which is impossible. Thus $O_3(G) = 1$.

Let $T$ be a Sylow 3-subgroup of $G$. Then $T \not\triangleleft G$ and so $j = |G : N_G(T)| \neq 1$. Since $j \equiv 1 \pmod{3}$ and $j | 8$, we have $j = 4$, whence $|N_G(T)| = 6$. The group $G$ acts on the right cosets of $N = N_G(T)$ by right multiplication. This determines a homomorphism $G \to S_4$ with kernel $K = \bigcap_{g \in G} N^g$. Since $T \not\triangleleft G$, $|T| = 3$ and $|N| = 6$, it is clear that $N \not\triangleleft G$. Therefore $K$ is proper in $N$, and so $|K|$ is 1 or 2 or 3. But $|K| = 3$ is impossible because $O_3(G) = 1$. If $|K| = 2$ then, since $K, T \triangleleft N$ and $K \cap T = 1$, we have $KT \simeq Z_6$, which is also impossible. Thus $|K| = 1$ and so $G \simeq S_4$. □

It follows that (2.11) is not valid in every group of order 24, except $S_4$. This completes the proof of Proposition 4. □

**Note.** It was claimed in [PK] that $\omega_{24}$, (2.10), (2.11), plus the formula

$$\big[\bigvee_{1 \leq i \leq 9} x_i^3 = 1\big] \vee \big[\bigvee_{1 \leq i < j \leq 9} (x_i x_j^{-1})^3 = 1\big]$$

form a weak basis of d-identities of $S_4$ (this result was also announced in [P], p.80). This is not quite true, for it is easy to see that the quaternion group $Q_8$ satisfies all these d-identities, has a faithful irreducible representation, but is not a section of $S_4$.

## 3. Disjunctive identities of $A_5$ and $S_5$

In the previous section we considered groups of very small orders, and it is not surprising that our arguments were so short and elementary. In this section we deal with groups of slightly higher order, and immediately the difficulty of the proofs increases.

In the sequel, we often use the following immediate consequence from the Sylow Theorem: if a set of d-identities eliminates all groups of order $p^k$ with $p$ a prime, then the same set eliminates all groups of the order divisible by $p^k$. Also, we sometimes informally say that the permutation $(123)(45)$ has a (cycle) form $(***)(**)$, the permutation $(15)(32)$ has form $(**)(**)$, etc.

THE GROUP $A_5$. 1. Every nonunit element of $A_5$ has one of the following three cycle forms: $(*****)$, $(***)$, $(**)(**)$. Therefore $A_5$ satisfies the d-identity

$$(x^2 = 1) \vee (x^3 = 1) \vee (x^5 = 1) \tag{3.1}$$

and, of course, $\omega_{60}$. These two d-identities may both be valid in the groups of the following orders only:

1, 2, 3, 4, 5, 6, 8, 9, 10, 12, 15, 16, 18, 20, 24, 25, 27, 30, 32, 36, 40, 45, 48, 50, 54, 60.

The groups of order 1,2,3,5 are all subgroups of $A_5$. The cyclic group of order 4 does not satisfy (3.1), but the noncyclic one is a subgroup of $A_5$. Similarly, the cyclic group of order 6 does not satisfy (3.1), but the noncyclic one (namely, $S_3$) is a subgroup of $A_5$: it can be generated, for example, by (123) and (12)(45). All groups of order 8, except $Z_2^3$, have elements of order 4 and do not satisfy (3.1).

2. The group $Z_2^3$ requires more efforts. First we note that if

$$a = (a_1 a_2)(a_3 a_4) \quad \text{and} \quad b = (b_1 b_2)(b_3 b_4)$$

are two elements of order 2 in $A_5$, then there are three possibilities:

(i) $\mathrm{Supp}(a) = \mathrm{Supp}(b)$, say $a = (12)(34)$ and $b = (13)(24)$. Then $ab = ba = (14)(23)$, so $a$ and $b$ are contained in the same Sylow 2-subgroup of $A_5$.

(ii) $|\mathrm{Supp}(a) \cap \mathrm{Supp}(b)| = 3$ and $a$ and $b$ have a common cycle, say $a = (12)(34)$ and $b = (12)(35)$. Then $ab = (345)$, so $|ab| = 3$.

(iii) $|\mathrm{Supp}(a) \cap \mathrm{Supp}(b)| = 3$ and $a$ and $b$ have no common cycles, say $a = (12)(34)$ and $b = (15)(24)$. Then $ab = (14325)$, so $|ab| = 5$.

Now consider the following formula in three variables $x_1, x_2, x_3$:

$$[\bigvee_i x_i^{15} = 1] \vee [\bigvee_{i<j} (x_i x_j)^{15} = 1] \vee [\bigvee_{(i,j,k) \neq 0} x_1^i x_2^j x_3^k = 1] \tag{3.2}$$

where the last disjunction is taken over all nonzero $\{0,1\}$-vectors $(i,j,k)$. Take any $x_1, x_2, x_3 \in A_5$. If the first two clauses of (3.2) are false then all $x_i$ are elements of order 2 belonging to the same Sylow 2-subgroup $H$ of $A_5$. Since $H$ is a 2-dimensional vector space over $GF(2)$, the last clause must hold. Thus (3.2) is valid in $A_5$. On the other hand, it is obviously false in $Z_2^3$.

Let $G$ be any group whose order is divisible by 8 and which satisfies (3.1). Then $G$ has a subgroup $H$ of order 8; in view of (3.1), $\exp H = 2$ and so $H \simeq Z_2^3$. This eliminates all groups of orders 16, 24, 32, 40 and 48.

3. There are two groups of order 9: $Z_9$ and $Z_3 \times Z_3$. The first does not satisfy (3.1). To eliminate the second, take the formula

$$(x_1^{10} = 1) \vee (x_2^{10} = 1) \vee (x_1 = x_2) \vee ((x_1 x_2)^{10} = 1) \vee ((x_1^2 x_2)^2 = 1). \tag{3.3}$$

If $x_1, x_2 \in A_5$ and $x_i^{10} \neq 1$, then $x_i$ are cycles of length 3. Without loss of generality we assume that $x_1 = (123)$. If $x_2$ moves the same symbols 1,2,3, then either $x_1 = x_2$ or $x_1 x_2 = 1$. If $x_2$ moves only one of the symbols 1,2,3, then $x_1 x_2$ is a cycle of length 5 and $(x_1 x_2)^5 = 1$. If $x_2$ moves two of the symbols 1,2,3, then there are two essentially different cases: $x_2 = (124)$ and $x_2 = (214)$. In the first case

$$x_1 x_2 = (123)(124) = (14)(23) \quad \text{and so} \quad (x_1 x_2)^2 = 1,$$

in the second case

$$x_1^2 x_2 = (132)(214) = (13)(24) \quad \text{and so} \quad (x_1^2 x_2)^2 = 1.$$

It follows that (3.3) is valid in $A_5$. Clearly it is not valid in $Z_3 \times Z_3$, and so the groups of order 9 are eliminated. The same is true for groups of orders 18, 27, 36, 45 and 54.

4. There are two groups of order 10: $Z_{10}$ and $D_{10}$. The first does not satisfy (3.1) but the second is a subgroup of $A_5$. There are five groups of order 12: $Z_{12}$, $Z_2 \times Z_2 \times Z_3$, $A_4$, $D_{12}$ and also the group $H = \langle a, b \,|\, a^3 = b^4 = 1,\ a^b = a^{-1}\rangle$. All of them except $A_4$ have elements of order 6 and therefore do not satisfy (3.1), but $A_4$ is a subgroup of $A_5$. The only group of order 15 is cyclic. Every group $G$ of order 20 has only one Sylow 5-subgroup, and so $G = Z_5 \leftthreetimes T$ where $|T| = 4$. Since $G$ may not have elements of order 4, we have $T = Z_2 \times Z_2$. But then $C_T(Z_5) \neq 1$, and so $G$ contains an element of order 10 which again contradicts (3.1). Thus the groups of order 10, 12, 15 and 20 are eliminated.

5. Let $G$ be a group of order 30. It is solvable and has a composition series

$$1 \triangleleft A \triangleleft B \triangleleft G$$

with factors of orders 2, 3 and 5. Therefore the order of $B$ is either 15, or 10, or 6. If $|B| = 15$, then $B = Z_{15}$ and (3.1) fails. If $|B| = 10$, then either $B \simeq Z_{10}$ or $B \simeq D_{10}$. The first case is again impossible, but in the second case $A \simeq Z_5$ is normal subgroup of $G$ and so $G/A$ is an extension of a group $B/A$ of order 2 by a group of order 3. But then $G/A \simeq Z_6$ and (3.1) fails.

Finally, if $|B| = 6$, then we may assume that $B \simeq S_3$. In this case $A \simeq A_3$ is characteristic in $B$ and so $A \triangleleft G$. Therefore $|G/A| = 15$, whence $G/A \simeq Z_{15}$. This eliminates all groups of order 30. Similar arguments apply to solvable groups of order 60; the only nonsolvable group of order 60 is $A_5$ itself.

6. To eliminate the groups of order 25, consider the formula

$$(x_1^6 = 1) \vee (x_2^6 = 1) \vee ((x_1 x_2)^6 = 1)$$
$$\vee ((x_1 x_2^2)^6 = 1) \vee ((x_1 x_2^3)^6 = 1) \vee ((x_1 x_2^4)^6 = 1). \quad (3.4)$$

It is valid in $A_5$. Indeed, take $x_1, x_2 \in A_5$ with $x_i^6 \neq 1$. Then the $x_i$ are 5-cycles. By Lemma 1, there exists $d$ ($1 \leq d \leq 4$) such that $x_1 x_2^d$ is not a 5-cycle, that is $(x_1 x_2^d)^6 = 1$. On the other hand, (3.4) fails in any group of order 25: take an element $x_2 \in G$ of order 5 and any element $x_1$ not in $\langle x_2 \rangle$. This eliminates the groups of order 25 and 50 and completes our argument.

**Theorem 1.** *Formulas $\omega_{60}$ and (3.1)–(3.4) form a basis of d-identities of $A_5$.* □

THE GROUP $S_5$. 1. We begin with two standard d-identities: $\omega_{120}$ and

$$(x^4 = 1) \vee (x^5 = 1) \vee (x^6 = 1). \quad (3.5)$$

They both are valid in $S_5$, and if they are valid in some group $G$, then the order of $G$ must be among the numbers

$$1, 2, 3, 4, 5, 6, 8, 9, 10, 12, 15, 16, 18, 20, 24, 25, 27, 30, 32, 36, 40, 45, 48,$$

$$50, 54, 60, 64, 72, 75, 80, 81, 90, 96, 100, 108, 120.$$

Noting that all groups of order $\leq 6$ are subgroups of $S_5$, consider the groups of order 8. Among them, $Z_8$ does not satisfy (3.5), but $D_8$ is a subgroup of $S_5$. To eliminate $Z_2 \times Z_4$ and $Q_8$, we take the formula

$$(x_1^{30} = 1) \vee (x_2^{30} = 1) \vee ((x_1 x_2)^{15} = 1) \vee ((x_1 x_2^2)^{15} = 1) \vee ((x_1 x_2^3)^{15} = 1). \quad (3.6)$$

It is easy to see that it is not valid in either $Z_2 \times Z_4$ or $Q_8$. We claim that it holds in $S_5$. Take $x_1, x_2 \in S_5$ such that $x_i^{30} \neq 1$. Then both $x_1$ and $x_2$ are 4-cycles, and the claim will follow from

**Lemma 3.** *If $x_1$ and $x_2$ are 4-cycles in $S_5$, then*

$$((x_1 x_2)^{15} = 1) \vee ((x_1 x_2^2)^{15} = 1) \vee ((x_1 x_2^3)^{15} = 1).$$

Proof. Without loss of generality we may assume that $x_1 = (1234)$, and let $x_2 = (a_1 a_2 a_3 a_4)$. Consider two cases.

(i) $\{a_1, a_2, a_3, a_4\} = \{1, 2, 3, 4\}$. We may assume that $x_2 = (1***)$, and then there are 6 choices for $x_2$. A direct verification shows that, unless $x_2 = x_1^{\pm 1}$, $x_1 x_2$ is a 3-cycle and then $(x_1 x_2)^{15} = 1$. If $x_2 = x_1^{-1}$ then $x_1 x_2 = 1$, but if $x_2 = x_1$ then $x_1 x_2^3 = x_1^4 = 1$.

(ii) $\{a_1, a_2, a_3, a_4\} \neq \{1, 2, 3, 4\}$. We may assume that $\{a_1, a_2, a_3, a_4\} = \{2, 3, 4, 5\}$ and $x_2 = (***5)$. Again there are 6 choices for $x_2$. For each of them, a direct verification shows that $x_1 x_2^k$ is a 5-cycle for some $k = 1, 2, 3$. □

2. The next step is to eliminate $Z_2^3$. As usual, it is the most painful. Recall that every Sylow 2-subgroup of $S_5$ is isomorphic to $D_8$ and therefore contains precisely two 4-cycles which uniquely determine it. Namely, the Sylow 2-subgroup containing a cycle $(abcd)$ consists of the elements

$$1, (abcd), (adcb), (ac), (bd), (ab)(cd), (ac)(bd), (ad)(bc).$$

We will denote it by Syl($abcd$). From this description, it is easy to deduce that if $y_1, y_2, y_3$ are 2-elements of $S_5$ such that any two of them are contained in one Sylow subgroup, then all three of them are contained in one Sylow subgroup.

Take three variables $x_1, x_2, x_3$ and denote $x_i^{15} = y_i$. Let $\theta(x_1, x_2, x_3)$ be the d-identity (2.4). Consider the formula

$$\theta(y_1, y_2, y_3) \vee \left[\bigvee_i y_i = 1\right] \vee \left[\bigvee_{i \neq j} (y_i y_j)^{15} = 1\right]$$

$$\vee \left[\bigvee_{j \neq i \neq k} (y_i y_j y_k)^{15} = 1\right] \vee \left[\bigvee_{j \neq i \neq k} (y_i y_j y_k y_j)^{15} = 1\right], \quad (3.7)$$

where $j \neq i \neq k$ means that $j \neq i$ and $i \neq k$ (but $j$ and $k$ may coincide). A basis of $Z_2^3$ refutes (3.7) in this group. Our goal is to show that (3.7) is valid in $S_5$. This is trivially true if $y_i = x_i^{15} = 1$ for some $i$, so we assume that $y_i \neq 1$ for all $i$. Hence the $y_i$ are elements of order 2 or 4.

Next, suppose that $y_1, y_2, y_3$ lie in the same Sylow 2-subgroup of $S_5$. Since this subgroup is isomorphic to $D_8$ and $\theta$ is a d-identity of $D_8$, it follows that $\theta(y_1, y_2, y_3)$ is valid, and so is (3.7). Thus we may assume that none of the Sylow 2-subgroups of $S_5$ contains all three $y_i$ simultaneously. It follows from the above property of the Sylow 2-subgroups that *two* of the elements $y_i$ do not belong to the same Sylow 2-subgroup. From now on, we assume that these elements are $y_1$ and $y_2$. There are three cases.

(i) *Both $y_1$ and $y_2$ are 4-cycles.* Then by Lemma 3 we have $((y_1y_2)^{15} = 1) \vee ((y_1y_2^2)^{15} = 1) \vee ((y_1y_2^3)^{15} = 1)$, and (3.7) holds.

(ii) *Only one of the $y_1$ and $y_2$ is a 4-cycle.* Let for example $y_1 = (1234)$. Then there are two possibilities.

a) $y_2$ is a transposition. Note that $y_2$ can not be equal to (13) or (24) because both these transpositions belong to $\mathrm{Syl}(1234)$. If $\mathrm{Supp}(y_2) \subset \mathrm{Supp}(y_1)$ then $y_1y_2$ is a 3-cycle (examples: $(1234)(12) = (234)$, $(1234)(14) = (123)$). If $\mathrm{Supp}(y_2) \not\subset \mathrm{Supp}(y_1)$ then, without loss of generality, we may assume that $y_2 = (15)$, whence $y_1y_2 = (12345)$. In all cases $(y_1y_2)^{15} = 1$ and (3.7) holds.

b) $y_2$ has form $(**)(**)$. Since $y_2 \notin \mathrm{Syl}(1234)$, it must involve symbol 5. Essentially, there may be only three possibilities:

$$y_2 = (15)(23) \quad \text{and then} \quad y_2y_1^2 = (15342),$$
$$y_2 = (15)(24) \quad \text{and then} \quad y_2y_1^2 = (153),$$
$$y_2 = (15)(34) \quad \text{and then} \quad y_2y_1^2 = (15324).$$

In all cases $(y_2y_1^2)^{15} = 1$ and (3.7) is valid.

(iii) *None of the $y_1$ and $y_2$ is a 4-cycle.* Then both $y_1$ and $y_2$ are involutions and therefore the subgroup they generate in $S_5$ is a dihedral group:

$$\langle y_1, y_2 \rangle = \langle y_1y_2 \rangle \rtimes \langle y_2 \rangle \simeq D_{2m},$$

where $\langle y_1y_2 \rangle \simeq Z_m$. Since $y_1$ and $y_2$ do not lie in the same Sylow 2-subgroup, $y_1y_2$ may not be a 2-element. It follows that the order of $y_1y_2$ is either 3, or 5, or 6. In the first two cases $(y_1y_2)^{15} = 1$. In the last case $y_1y_2 = (***)(**)$, from which one can easily see that

$$y_1 = (**)(**), \quad y_2 = (**), \quad \mathrm{Supp}(y_1) \cup \mathrm{Supp}(y_2) = \{1,2,3,4,5\}.$$

Without loss of generality we may assume that $y_1 = (12)(34)$ and $y_2 = (15)$. Again we have two possibilities.

a) $y_3$ is a transposition. If $y_2$ and $y_3$ overlap, then $(y_2y_3)^3 = 1$. Otherwise there are two essentially different cases:

$$y_3 = (23) \quad \text{and then} \quad y_1y_3y_2 = (13425),$$
$$y_3 = (34) \quad \text{and then} \quad y_1y_3y_2 = (125).$$

In both cases $(y_1 y_3 y_2)^{15} = 1$.

b) $y_3$ has form $(**)(**)$. Suppose first that $\text{Supp}(y_3) = \{1,2,3,4\}$. If $y_3 = y_1$, then $y_1 y_3 = 1$ and there is nothing to prove. All other possibilities are essentially equivalent to the case when $y_3 = (13)(24)$. Then we have $y_1 y_2 y_3 y_2 = (14532)$ and so $(y_1 y_2 y_3 y_2)^{15} = 1$. Now suppose that $\text{Supp}(y_3) \neq \{1,2,3,4\}$. Then we have two essentially different cases ($y_2$ can now be disregarded):

$$y_3 = (15)(34) \quad \text{and then} \quad y_1 y_3 = (125),$$
$$y_3 = (15)(24) \quad \text{and then} \quad y_1 y_3 = (14325).$$

In both cases $(y_1 y_3)^{15} = 1$. Thus we have proved that (3.7) is a d-identity of $S_5$. This concludes the consideration of groups of order 8.

3. Consider the d-identity

$$(x_1^{20} = 1) \vee (x_2^{20} = 1) \vee (x_1^2 = x_2^2) \vee ((x_1^2 x_2^2)^5 = 1) \vee ((x_1^4 x_2^2)^2 = 1). \tag{3.8}$$

Note that if $\delta(x_1, x_2)$ is the formula (3.3), then $\delta(x_1^2, x_2^2)$ is precisely (3.8). Since (3.3) is a d-identity of $A_5$ and $x_i^2 \in A_5$ for any $x_i \in S_5$, we see that (3.8) is valid in $S_5$. On the other hand, it obviously fails in the groups of order 9 and so in the groups of orders 18, 27, 36, 45, 54, 72, 81, 90, 108.

4. Among the two groups of order 10, the cyclic one does not satisfy (3.5), but the noncyclic one is a subgroup of $S_5$. Among the groups of order 12, $Z_{12}$ does not satisfy (3.5), but $A_4$ and $D_{12}$ are subgroups of $S_5$. To eliminate $Z_2 \times Z_6$, consider the d-identity

$$(x_1^{15} = 1) \vee (x_2^{15} = 1) \vee (x_1^4 = 1) \vee (x_2^4 = 1) \vee ((x_1 x_2)^{15} = 1) \vee ((x_1 x_2)^4 = 1). \tag{3.9}$$

Let $x_1, x_2 \in S_5$ and suppose that the first four clauses of (3.9) are false. Then both $x_1$ and $x_2$ are of order 6, i.e. have form $(***)(**)$. But then $x_1 x_2$ belongs to $A_5$, and so either $(x_1 x_2)^4 = 1$ or $(x_1 x_2)^{15} = 1$. Thus (3.9) holds in $S_5$. On the other hand, it fails in $Z_2 \times Z_6$: if $x_1 = (1, b)$ and $x_2 = (a, b)$, where $a, b$ is a natural basis of $Z_2 \times Z_6$, then $x_1, x_2, x_1 x_2$ are all elements of order 6.

The remaining group of order 12 is

$$G = \langle a, b \mid a^3 = b^4 = 1, a^b = a^{-1} \rangle = Z_3 \rtimes Z_4.$$

It is eliminated by the following d-identity from [PK] (see also [P], p.81):

$$(x_1^{10} = 1) \vee (x_2^{10} = 1) \vee (x_1^4 = x_2^4) \vee (x_1^6 = x_2^6)$$
$$\vee ((x_1^2 x_2^2)^3 = 1) \vee ((x_1^2 x_2^2)^4 = 1) \vee ((x_1^2 x_2^2)^5 = 1). \tag{3.10}$$

Leting $x_1 = a$, $x_2 = b$, we see that it fails in $G$. It is proved in [PK] that (3.10) holds in $S_5$ (if the first four clauses of (3.10) fail, then one of the $x_i$ has order 3 or 6 and the other has order 4; all such possibilities can be checked one by one).

5. Now we consider groups of orders 15, 30, 60, 120. The only group of order 15 is $Z_{15}$ which does not satisfy (3.5). Every solvable group of order $15m$ with $(15, m) = 1$ has a subgroup of order 15, and therefore all solvable groups of orders 30, 60, 120 are also eliminated. The only nonsolvable group of order $\leq 60$ is $A_5$ which is a subgroup of $S_5$. Finally, there are only three nonsolvable groups of order 120, namely $S_5$, $SL(2, 5)$ and $A_5 \times Z_2$, but the latter two have elements of order 10 and do not satisfy (3.5).

6. Every group of order 16 has an abelian section $A$ of order 8. If $\exp A \leq 4$, then $A$ is either $Z_2 \times Z_4$ or $Z_2^3$ and hence does not satisfy (3.6) or (3.7), respectively. Otherwise $A$ does not satisfy (3.5). This eliminates all groups of orders 16, 32, 48, 64, 80, 96.

7. Every group $G$ of order 20 has a unique Sylow 5-subgroup. It follows that $G = H \rtimes T$, where $H \simeq Z_5$, $|T| = 4$. Let $C = C_T(H)$. If $C \neq 1$, then $G$ has an element of order 10 and does not satisfy (3.5). If $C = 1$, then it is easy to see that $T \simeq Z_4$ and $G$ is isomorphic to the subgroup of $S_5$ generated by (12345) and (2354).

Similarly, every group of $G$ order 40 can be presented as $G = H \rtimes T$, where $H \simeq Z_5$, $|T| = 8$. Since $\operatorname{Aut} Z_5 \simeq Z_4$, we see that $C_T(H) \neq 1$ and so $G$ must have an element of order 10.

8. The next step is to eliminate groups of orders divisible by 25. This can be done by the d-identity

$$(x_1^{12} = 1) \vee (x_2^{12} = 1) \vee ((x_1 x_2)^{12} = 1)$$
$$\vee ((x_1 x_2^2)^{12} = 1) \vee ((x_1 x_2^3)^{12} = 1) \vee ((x_1 x_2^4)^{12} = 1), \quad (3.11)$$

using the same arguments as for the group $A_5$ (cf. (3.4)). Moreover, it is easy to see that (3.11) is valid *even in* $S_6$. Indeed, if $x_1, x_2 \in S_6$ and $x_i^{12} \neq 1$, then $x_i$ are 5-cycles, and the rest follows from Lemma 1.

9. It remains to consider the groups of order 24. We will prove that every group $G$ of order 24 satisfying the above d-identities is a section of $S_5$. We may assume that $G$ has no elements of order 12 (by (3.5)) and that all of its sections of order 8 are isomorphic to $D_8$ (because the other groups of order 8 have been eliminated).

The proof is similar to that of Lemma 2. Note that the number of Sylow 3-subgroups of $G$ may be 1 or 4. In the first case

$$G = A \rtimes B \quad \text{where} \quad A \simeq Z_3, \ B \simeq D_8.$$

Let $C = C_B(A)$, then $|B/C|$ is 1 or 2. If $|B/C| = 1$, then $G = A \times B$ and $G$ must have an element of order 12, which is impossible. If $|B/C| = 2$, then $|C| = 4$ and $\langle A, C \rangle = A \times C$. If $C \simeq Z_4$, we have an element of order 12. If $C \simeq Z_2 \times Z_2$, then $\langle A, B \rangle \simeq Z_3 \times Z_2 \times Z_2 = Z_6 \times Z_2$, which was eliminated earlier.

Suppose now that the number of Sylow 3-subgroups of $G$ is 4. Let $T$ be one of these subgroups and let $N = N_G(T)$, $K = \cap_{g \in G} N^g$. Repeating the proof of Lemma 2, we eliminate all cases except $|K| = 2$. Then $|G/K| = 12$ and, since $G/K \subseteq S_4$ and the only subgroup of order 12 of $S_4$ is $A_4$, we obtain $G/K \simeq A_4$. Let $B/K$ be the unique Sylow

2-subgroup of $G/K$. Then $|B| = 8$ and $G = B \rtimes H$ with $|H| = 3$. By initial assumption, $B \simeq D_8$. Therefore the group $H/C_H(B)$ embeds into $\operatorname{Aut} D_8$, which is a 2-group. It follows that $H/C_H(B) = 1$, whence $G = B \times H \simeq D_8 \times Z_3$, and we again have an element of order 12, which is impossible. This completes the proof.

**Theorem 2.** *Formulas $\omega_{120}$ and (3.5)–(3.11) form a basis of d-identities of $S_5$.* $\square$

**Note.** The first study of d-identities of $S_5$ was done by Plotkin and Kushkuley [PK] (see also [P; Section 2.5.1]). They exibited a list of seven formulas and claimed that it was a weak basis of d-identities of $S_5$. This result contains useful information and has partially influenced our arguments. One should note, however, that it has two gaps. First, the list provided in [PK] and [P] does not eliminate the quaternion group $Q_8$, which is an irreducible linear group, but is not a section of $S_5$. Second, one of the formulas in that list, namely

$$(x_1^3 = 1) \vee (x_2^3 = 1) \vee (x_1^5 = 1) \vee (x_2^5 = 1) \vee (x_1^4 = x_2^4)$$
$$\vee ((x_1 x_2)^3 = 1) \vee ((x_1 x_2)^4 = 1) \vee ((x_1 x_2)^5 = 1) \quad (3.12)$$

is not valid in $S_5$ at all: take, for example, $x_1 = (123)(45)$, $x_2 = (14)(25)$.

## 4. Disjunctive identities of $A_6$

The procedure of computing a basis of d-identities of $A_6$ is more comlicated. If we did it "from scratch", it would be substantially longer. Fortunately, a number of intermediate steps have been already done in the previous sections. We begin with stating the final result.

**Theorem 3.** *The following formulas form a basis of d-identities of $A_6$:*

$$\omega_{360} = \bigvee_{0 \leq i < j \leq 360} (x_i = x_j), \quad (4.1)$$

$$(x^3 = 1) \vee (x^4 = 1) \vee (x^5 = 1), \quad (4.2)$$

$$\big[\bigvee_{1 \leq i \leq 9} x_i^{20} = 1\big] \vee \big[\bigvee_{i<j} x_i = x_j\big] \vee \big[\bigvee_{i<j}(x_i x_j)^{20} = 1\big] \vee \big[\bigvee_{i<j}(x_i x_j^2)^{20} = 1\big], \quad (4.3)$$

$$(x_1^{30} = 1) \vee (x_2^{30} = 1) \vee \left((x_1 x_2)^{15} = 1\right) \vee \left((x_1 x_2^2)^{15} = 1\right) \vee \left((x_1 x_2^3)^{15} = 1\right), \quad (4.4)$$

$$\theta(x_1, x_2, x_3) \vee \big[\bigvee_i x_i^{15} = 1\big] \vee \big[\bigvee_{i \neq j}(x_i x_j)^{15} = 1\big]$$
$$\vee \big[\bigvee_{i \neq j}(x_i x_j^2)^{15} = 1\big] \vee \big[\bigvee_{i \neq j}(x_i x_j^3)^{15} = 1\big], \quad (4.5)$$

$$(x_1^{12} = 1) \vee (x_2^{30} = 1) \vee ((x_1 x_1^{x_2})^{12} = 1)$$
$$\vee ((x_1 x_1^{3x_2})^{12} = 1) \vee ((x_1^2, x_1^{x_2})^{12} = 1), \quad (4.6)$$

$$(x_1^{12} = 1) \vee (x_2^{12} = 1) \vee ((x_1 x_2)^{12} = 1)$$
$$\vee ((x_1 x_2^2)^{12} = 1) \vee ((x_1 x_2^3)^{12} = 1) \vee ((x_1 x_2^4)^{12} = 1). \quad (4.7)$$

**Notes.** 1. The d-identities (4.4) and (4.7) were introduced in the previous section (as (3.6) and (3.11)) and are repeated here for convenience.

2. As usual, $\theta(x_1, x_2, x_3)$ is the formula (2.4).

3. In (4.3), the indices $i, j$ run over the numbers 1,2,...,9; in (4.5) they run over the numbers 1,2,3.

P r o o f. A. First we prove that the d-identities (4.1) – (4.6) are valid in $A_6$. This is obvious for (4.1) and (4.2).

A1. To show that (4.3) holds in $A_6$, we first prove an auxiliary fact similar to Lemma 1.

**Lemma 4.** *Let $\alpha$ and $\beta$ be elements of order 3 in $S_6$ not belonging to the same Sylow 3-subgroup. Then either $\alpha\beta$ or $\alpha\beta^2$ is not an element of order 3.*

P r o o f is straightforward. Each Sylow 3-subgroup of $A_6$ is generated by two disjoint 3-cycles. For example, the Sylow 3-subgroup generated by (123) and (456) contains the following nonunit elements:

$$(123), \quad (132), \quad (456), \quad (465),$$

$$(123)(456), \quad (123)(465), \quad (132)(456), \quad (132)(465).$$

It follows that every Sylow 3-subgroup of $A_6$ is uniquely determined by any of its nonunit elements.

Let $\alpha$ and $\beta$ be 3-cycles. Then $s = |\mathrm{Supp}(\alpha) \cap \mathrm{Supp}(\beta)| = 1$ or 2, because otherwise $\alpha$ and $\beta$ lie in the same Sylow 3-subgroup. If $s = 1$, say $\alpha = (123)$ and $\beta = (145)$, then $\alpha\beta = (12345)$ is of order 5. If $s = 2$, then we may assume that $\alpha = (123)$ and $\mathrm{Supp}(\beta) = \{1, 2, 4\}$. Then there are only two choices for $\beta$ and for each of them either $\alpha\beta$ or $\alpha\beta^2$ is (14)(23), i.e. is of order 2.

Now let both $\alpha$ and $\beta$ be of the form $(***)(***)$. Let $\alpha = (123)(456)$ and $\beta = (b_1 b_2 b_3)(b_4 b_5 b_6)$. Since $\alpha$ and $\beta$ do not lie in one Sylow 3-subgroup, the set $\{b_1, b_2, b_3\}$ is neither $\{1, 2, 3\}$ nor $\{4, 5, 6\}$. Without loss of generality, we may assume that $\{b_1, b_2, b_3\} = \{1, 2, 4\}$. Then there are four choices for $\beta$:

$$\beta_1 = (124)(356), \qquad \beta_2 = (124)(365),$$
$$\beta_3 = (142)(356), \qquad \beta_4 = (142)(365).$$

Since $\beta_1^2 = \beta_4$ and $\beta_2^2 = \beta_3$, it remains to notice that $\alpha\beta_2 = (14326)$ and $\alpha\beta_4 = (26)(34)$.

In the case when only one of the $\alpha, \beta$ is a 3-cycle, the proof is similar and we omit it. □

Now we can prove that (4.3) is a d-identity of $A_6$. Let $x_1, x_2, \ldots, x_9 \in A_6$. If the first clause of (4.3) fails, then all $x_i$ are of order 3. If, in addition, the second clause fails, then they are pairwise different. Among nine different elements of order 3 in $A_6$, there must be two, $x_i$ and $x_j$, that do not belong to the same Sylow 3-subgroup. But then, by Lemma 4, either the third or the fourth clause of (4.3) is valid.

A2. It was proved in the previous section that (4.4) (that is, (3.6)) is valid in $S_5$. Now we show that it is valid in $A_6$. If $x_1, x_2 \in A_6$ and $x_i^{30} \neq 1$, then both $x_1$ and $x_2$ have form $(***)(**)$. Let $x_1 = (1234)(56)$ and $x_2 = (a_1 a_2 a_3 a_4)(a_5 a_6)$. There are three possibilities.

(i) $\{a_1, a_2, a_3, a_4\} = \{1, 2, 3, 4\}$. In this case the proof coincides with that for the group $S_5$.

(ii) $\{a_1, a_2, a_3, a_4\}$ contains exactly three of the numbers 1,2,3,4. Suppose these numbers are 2,3,4 and $x_2 = (***5)(16)$. Then there are six choices for $x_2$, and for each of them one can directly verify that

$$(x_1 x_2^k)^{15} = 1 \quad \text{for some} \quad k \leq 3.$$

(iii) $\{a_1, a_2, a_3, a_4\}$ contains exactly two of the numbers 1,2,3,4. Then there are six essentially different choices for $x_2$:

$$(5612)(34), \quad (5613)(24), \quad (5621)(34),$$
$$(5162)(34), \quad (5163)(34), \quad (5261)(34).$$

Again a direct verification shows that for each of these choices $(x_1 x_2^k)^{15} = 1$ for some $k \leq 3$.

A3. It will be more difficult to prove the validity of (4.5). First we describe the structure of Sylow 2-subgroups of $A_6$. Just as in $S_5$, every such subgroup is isomorphic to $D_8$ and is uniquely determined by each of the elements of order 4 it contains. Namely, the Sylow 2-subgroup containing the permutation $\alpha = (abcd)(ef)$ consists of the elements

$$1, \quad (abcd)(ef), \quad (adcb)(ef),$$
$$(ac)(ef), \quad (bd)(ef), \quad (ab)(cd), \quad (ac)(bd), \quad (ad)(bc).$$

This subgroup will be denoted by $\mathrm{Syl}(\alpha)$. It is easy to see that the Sylow 2-subgroups of $A_6$ possess the same property as those of $S_5$: if $x_1, x_2, x_3$ are 2-elements of $A_c$ and any two of them lie in one Sylow 2-subgroup, then all three of them lie in one Sylow 2-subgroup.

Choose $x_1, x_2, x_3 \in A_6$ such that $x_i^{15} \neq 1$. Then the $x_i$ are 2-elements. If they lie in one Sylow 2-subgroup, then $\theta(x_1, x_2, x_3)$, being a d-identity of $D_8$, must be valid. Otherwise two of the $x_i$, say $x_1$ and $x_2$, do not belong to the same Sylow 2-subgroup. Then there are three possibilities.

(i) *Both $x_1$ and $x_2$ are of order 4.* We may assume that $x_1 = (1234)(56)$ and $x_2 = (abcd)(ef)$. Consider three cases.

  a) $\{a,b,c,d\} = \{1,2,3,4\}$. In this case one can repeat the proof of Lemma 3, part (i), and show that

$$\left((x_1, x_2)^{15} = 1\right) \vee \left((x_1 x_2^2)^{15} = 1\right) \vee \left((x_1 x_2^3)^{15} = 1\right). \tag{4.8}$$

  b) $\{a,b,c,d\}$ contains exactly one of the numbers 5,6. Assuming that $x_2 = (5***)(64)$, we have six choices for $x_2$, and for each of them one can directly verify that (4.8) is valid.

  c) $\{a,b,c,d\}$ contains both 5 and 6. Without loss of generality, we may assume that either $x_2 = (516*)(**)$ or $x_2 = (561*)(**)$. This leads to six different choices for $x_2$

$$(5612)(34), \quad (5613)(24), \quad (5614)(23),$$
$$(5162)(34), \quad (5163)(24), \quad (5164)(23),$$

and again a direct verification shows that (4.8) is valid.

(ii) *Only one of the $x_1$ and $x_2$ is of order 4.* Assume that $x_1 = (1234)(56)$. Then $x_2$ has form $(**)(**)$ and, since $x_2 \notin \mathrm{Syl}(x_1)$, it follows that $\mathrm{Supp}(x_2) \neq \{1,2,3,4\}$. Hence $\mathrm{Supp}(x_2)$ contains at least one of the numbers 5 and 6.

  a) Let $\mathrm{Supp}(x_2)$ contain precisely one of these numbers, say 5. Repeating the proof of Theorem 2, Step 2, (ii), b), we obtain that $(x_2, x_1^2)^{15} = 1$, whence (4.5) holds.

  b) Now suppose that $x_2 = (56)(**)$. Since $(13)(56)$ and $(24)(56)$ lie in $\mathrm{Syl}(x_1)$, there are four choices for $x_2$.

$$(56)(12), \quad (56)(14), \quad (56)(13), \quad (56)(34).$$

In each case $x_1 x_2$ is a 3-cycle, so $(x_1 x_2)^3 = 1$ and (4.5) holds.

  c) Finally, let $x_2 = (5*)(6*)$. There are only two essentially different choices for $x_2$: $(51)(62)$ and $(51)(63)$. In the first case, $x_2 x_1^2 = (153)(264)$, in the second case $x_2 x_1 = (164)(235)$, so (4.5) is valid again.

(iii) *None of the $x_1$ and $x_2$ is of order 4.* Then they both are involutions and therefore $\langle x_1, x_2 \rangle$ is isomorphic to a dihedral group:

$$\langle x_1, x_2 \rangle = \langle x_1 x_2 \rangle \rtimes \langle x_2 \rangle \cong D_{2m}$$

where $\langle x_1 x_2 \rangle \cong Z_m$. Since $x_1$ and $x_2$ are not contained in the same Sylow 2-subgroup, $x_1 x_2$ is not a 2-element. Therefore $(x_1 x_2)^{15} = 1$ and (4.5) holds.

A4. To show that (4.6) holds in $A_6$, take $x_1, x_2 \in A_6$ such that $x_1^{12} \neq 1$ and $x_2^{30} \neq 1$. Then $x_1$ is a 5-cycle and $x_2$ has form $(***)(**)$. Without loss of generality, we may assume that $x_1 = (12345)$, then $x_1^{x_2}$ is a 5-cycle involving 6. The rest follows from

**Lemma 5.** *Let $\alpha$ and $\beta$ be 5-cycles in $S_6$ such that $Supp(\alpha) \neq Supp(\beta)$. Then one of the permutations $\alpha\beta$, $\alpha\beta^3$, $\alpha^2\beta$ is not a 5-cycle.*

Proof. Unfortunately, we could not find a nice proof of this fact (which certainly exists), and simply considered all possible cases. First, we may assume that $\alpha = (12345)$ and $\beta = (b_1 b_2 b_3 b_4 6)$, where $b_i \neq 5$. Thus there are 24 choices for $(b_1 b_2 b_3 b_4)$. Next, we notice that

$$\text{if } \beta = (1**26) \quad \text{then} \quad \alpha\beta = (16)(****),$$
$$\text{if } \beta = (2**36) \quad \text{then} \quad \alpha\beta = (26)(****),$$
$$\text{if } \beta = (3**46) \quad \text{then} \quad \alpha\beta = (36)(****).$$

It remains to consider 18 possibilities for $(b_1 b_2 b_3 b_4)$, which can be done by direct verification. For example:

$$\text{if } \beta = (12346) \quad \text{then} \quad \alpha\beta = (136)(245),$$
$$\text{if } \beta = (13246) \quad \text{then} \quad \alpha\beta^3 = (26)(45),$$
$$\text{if } \beta = (21346) \quad \text{then} \quad \alpha^2\beta = (1435)(26),$$

and so on. $\square$

Lemma 5 implies that (4.6) is a d-identity of $A_6$. From the previous section, we know that (4.7) is also a d-identity of $A_6$. This completes the first part of the proof of Theorem 3.

B. Now we have to show that if a group $G$ satisfies the $d$-identities (4.1)–(4.6), then it is a section of $A_6$. Let $|G| = n$; it follows from (4.1) and (4.2) that

$$n \leq 360 \quad \text{and} \quad n = 2^k 3^\ell 5^m.$$

B1. We show that in fact $n$ must be a divisor of 360, that is, $k \leq 3$, $\ell \leq 2$, $m \leq 1$. Since the d-identity (4.7) is not valid in any group of order 25, it follows that $m \leq 1$. To show that $\ell \leq 2$ assume that a group $A$ of order 27 satisfies the d-identities (4.1)–(4.7). Then (4.2) implies that $\exp A = 3$. We can choose in $A$ nonunit elements $x_1, x_2, \ldots, x_9$ such that $x_i \neq x_j^{\pm 1}$ for $i \neq j$. Then $x_i x_j \neq 1$ and $x_i x_j^2 \neq 1$ (since $\exp A = 3$, we have $x_j^2 = x_j^{-1}$), and so (4.3) fails. It follows that $\ell \leq 2$.

To show that $k \leq 3$, consider first the abelian groups of order 8. Among them, $D_8$ is a subgroup of $A_6$, $Z_8$ does not satisfy (4.2), $Z_4 \times Z_2$ and $Q_8$ does not satisfy (4.4) (see the proof of Theorem 2, Step 1), and $Z_2^3$ does not satisfy (4.5) (take a basis of $Z_2^3$ for the $x_1, x_2, x_3$). Repeating the proof of Theorem 2, Step 6, we eliminate all groups of order 16. It follows that $k \leq 3$.

Thus, if a group $G$ of order $n$ satisfies the d-identities (4.1)–(4.7), then

$$n = 2^k 3^\ell 5^m, \quad \text{where} \quad k \leq 3, \quad \ell \leq 2, \quad m \leq 1,$$

that is, $n$ is one of the numbers

$$1,\ 2,\ 3,\ 4,\ 5,\ 6,\ 8,\ 9,\ 10,\ 12,\ 15,\ 18,\ 20,\ 24,$$

$$30,\ 36,\ 40,\ 45,\ 60,\ 72,\ 90,\ 120,\ 180,\ 360.$$

B2. All groups of order $\leq 6$, except $Z_6$, are subgroups of $A_6$, but $Z_6$ does not satisfy (4.2). The groups of order 8 have been already taken care of. The cyclic groups of order 9 and 10 do not satisfy (4.2), but the noncyclic ones are subgroups of $A_6$. All groups of order 12, except $A_4$, have elements of order 6 and do not satisfy (4.2), but $A_4$ is a subgroup of $A_6$. Repeating the proof of Theorem 2, Step 5, and noting that $S_5$ does not satisfy (4.2), we eliminate all groups of orders 15, 30, 60, and 120. Lemma 2 and (4.2) eliminate all groups of order 24 except $S_4$, but $S_4$ is a subgroup of $A_6$.

B3. *Groups of order 18.* Every group $G$ of order 18 splits into a semidirect product

$$G = A \rightthreetimes B, \quad \text{where} \quad |A| = 9 \quad \text{and} \quad |B| = 2.$$

If $A \simeq Z_9$, then (4.2) fails, so we assume that $A \simeq Z_3 \times Z_3$. Let $1 \neq b \in B$. If $a^b \neq a^{-1}$ for some $a \neq a \in A$, then $c = aa^b \neq 1$ and $c^b = a^b a = aa^b = c$. Hence $cb$ is an element of order 6 and so $G$ is eliminated by (4.2). If $a^b = a^{-1}$ for all $a \in A$, then $G$ is isomorphic to the subgroup of $A_6$ generated by (123), (456) and (12)(45).

B4. *Groups of order 36.* Every group $G$ of order 36 is solvable and therefore if $H$ is a maximal normal subgroup of $G$, then $|G : H| = 2$ or $|G : H| = 3$.

(i) Let $|G : H| = 3$. Since all groups of order 12 except $A_4$ have elements of order 6, we have $H \simeq A_4$. Hence

$$G = A \rightthreetimes T, \quad \text{where} \quad A \cong Z_2 \times Z_2 \quad \text{and} \quad |T| = 9.$$

Since $|\mathrm{Aut}(A)| = 6$, we see that $C_T(A)$ is nontrivial. It follows that $G$ must have elements of order 6, which is impossible.

(i) Let $|G : H| = 2$. Then $|H| = 18$ and, as established in B3, $H = A \rightthreetimes B$ with $A \simeq Z_3 \times Z_3$, $B \simeq Z_2$, and $a^b = a^{-1}$ for all $a \in A$. Therefore $A$ must be normal in $G$, whence $G = A \rightthreetimes T$, $|T| = 4$. If $C_T(A) \neq 1$, then $G$ again has elements of order 6. Hence $C_T(A) = 1$.

a) If $T \simeq Z_2 \times Z_2$, then

$$G = A \rightthreetimes (\langle b_1 \rangle \times \langle b_2 \rangle), \quad \text{where} \quad |b_i| = 2.$$

As shown in B3, both $b_i$ must act on $A$ as follows: $a^{b_i} = a^{-1}$ for all $a \in A$. Then $b_1 b_2^{-1}$ centralizes $A$, which is impossible.

b) If $T \simeq Z_4$, then $G = A \rightthreetimes \langle b \rangle$, where $|b| = 4$. Since $b^2$ is an involution, the same argument as in B3 shows that

$$\forall a \in A: \quad a^{b_2} = a^{-1}. \tag{4.9}$$

Take any $1 \neq a_1 \in A$. If $a_1^b = a_1$, we get an element of order 12. If $a_1^b = a_1^{-1}$, then $a_1^{b^2} = a_1$, contradicting (4.9). It follows that $a_1^b \neq a_1^{\pm 1}$, and then $a_1$ and $a_1^b$ form a basis of $A = Z_3 \times Z_3$. Denoting $a_1^b = a_2$, we have

$$G = (\langle a_1 \rangle \times \langle a_2 \rangle) \rtimes \langle b \rangle,$$

where $|a_i| = 3$, $|b| = 4$, $a_1^b = a_2$, $a_2^b = a_1^{b^2} = a_1^{-1}$. But this group is isomorphic to the subgroup of $A_6$ generated by $a_1 = (123)$, $a_2 = (456)$ and $b = (1425)(36)$, for

$$a_1^b = (5241)\,(36)\,(123)\,(1425)\,(36) = (456) = a_2,$$
$$a_2^b = (5241)\,(36)\,(456)\,(1425)\,(36) = (132) = a_1^{-1}.$$

This completes the consideration of groups of order 36.

B5. *Groups of orders 45, 90 and 180.* It follows directly from the Sylow Theorem that any group $G$ of order 45 is a direct product $G = A \times B$, where $|A| = 9$ and $|B| = 5$. Hence $G$ contains an element of order 15 and does not satisfy (4.2). Every group of order 90 is solvable. Since $90 = 45 \cdot 2$ and $(45, 2)=1$, this group has a subgroup of order 45 and therefore does not satisfy (4.2) either. The same argument works for solvable groups of order 180. There is only one unsolvable group of order 180, namely $A_5 \times Z_3$, which also contains elements of "wrong" orders.

B6. *Groups of orders 20 and 40.* Every group $G$ of order 20 can be presented as

$$G = A \rtimes B, \quad \text{where} \quad A = Z_5, \ |B| = 4.$$

If $C_B(A) \neq 1$, we get an element of order 10, which is impossible. Otherwise $B \simeq \mathrm{Aut}(Z_5)$ and $G$ can be defined by the presentation

$$G = \langle a, b \mid a^5 = b^4 = 1, \ a^b = a^2 \rangle.$$

But this group does not satisfy the d-identity (4.6): set $x_1 = a$ and $x_2 = b$. As to the groups of order 40, they are eliminated literally as in Theorem 2, Step 7.

B7. *Groups of order 72.* Let $G$ be a group of order 72 satisfying (4.1)–(4.7). Then $G$ is solvable and if $N$ is a maximal normal subgroup of $G$, then $|G:N| = 2$ or $|G:N| = 3$.

(i) Let $|G:N| = 2$. Then $|N| = 36$ and, as established in B4,

$$N = A \rtimes T, \quad \text{where} \quad A = Z_3 \times Z_3, \quad T = Z_4, \quad C_T(A) = 1$$

(recall that all other cases led to elements of order 6). Hence we have $A \triangleleft G$ and $|G/A| = 8$. All groups of order 8, except $D_8$, were eliminated by (4.2), (4.4) and (4.5). So $G/A \simeq D_8$ and, since $(|A|, |G/A|) = 1$, we have

$$G = A \rtimes D_8 = (Z_3 \times Z_3) \rtimes D_8.$$

But then $G$ has as a subgroup $(Z_3 \times Z_3) \rightthreetimes (Z_2 \times Z_2)$. Repeating B4 (ii), we get an element of order 6, which is impossible.

(ii) Let $|G : N| = 3$. Then $|N| = 24$ and, by Lemma 2, $N \simeq S_4$. Let $b \in G/N$. Clearly $|b|$ is divisible by 3, and it follows from (4.2) that $|b| = 3$. Therefore

$$G = N \rightthreetimes \langle b \rangle \simeq S_4 \rightthreetimes \langle b \rangle.$$

If $b$ centralizes $N$, we again have an element of order 6, otherwise $\langle b \rangle \subset \mathrm{Aut}(S_4)$. Since all automorphisms of $S_4$ are inner, it follows that $b$ acts on $S_4$ as conjugation by some 3-cycle, say (123):

$$\forall \alpha \in S_4 : \quad \alpha^b = (321)\,\alpha\,(123).$$

But then one can verify that the element $g = (12) \cdot b$ of the group $G = S_4 \rightthreetimes \langle b \rangle$ has order 6.

B8. *Groups of order 360.* Every solvable group of order $360 = 72 \cdot 5$ has a subgroup $H$ of order 72. It has just been proven that if $H$ satisfies (4.1)–(4.7), it must have elements of order 6.

Let $G$ be a nonsolvable group of order 360. If $G$ is simple, it is isomorphic to $A_6$. Otherwise it has composition series

$$1 \triangleleft A \triangleleft B \triangleleft G \tag{4.10}$$

with factors $Z_2, Z_3$ and $A_5$. If $|G/B| = 2$, then $B$ is a nonsolvable group of order 180. We already know that there is only one such group $A_5 \times Z_3$, which has elements of order 15 and 12. If $|G/B| = 3$, then $B$ is a nonsolvable group of order 120, i.e. one of the groups $S_5$, $\mathrm{SL}(2,5)$, $A_5 \times Z_2$. Each of them has elements of order $\geq 6$ and is excluded by (4.2).

Finally, let $G/B \simeq A_5$. Then $|B| = 6$ and, since elements of order 6 are prohibited, we have $B \simeq S_3$. If $C = C_G(B)$, then $G/C$ acts faithfully on $B$. Since $\mathrm{Aut}(S_3) \simeq S_3$, we get $|G/C| \leq G$. Thus $C$ must have elements of order 5, whence $G$ must have elements of order 10 and 15. This concludes the proof. $\square$

## 5. Disjunctive identities of the dihedral groups

Let $D_{2m} = \langle a, b \,|\, a^m = b^2 = 1, \ a^b = a^{-1} \rangle$ be the dihedral group of order $n = 2m$. In this section we find bases of $d$-identities of the groups $D_{2m}$ for any $m$. Since we have to deal with an infinite family of groups, rather than a concrete single group, the structure of the proof is different from that in the previous sections.

1. We begin with the following observation.

**Lemma 5.** *A group $G$ satisfies the $d$-identity*

$$(x^2 = 1) \vee (y^2 = 1) \vee (xy = yx) \tag{5.1}$$

*if and only if $G = A \rtimes B$, where $A$ is abelian and $B$ is either trivial or cyclic of order 2 acting by inversion on $A$.*

Proof. Clearly if $G$ has the described structure, it satisfies (5.1). On the other hand, let $G$ satisfy (5.1). If $\exp G = 2$, there is nothing to prove. Otherwise, let $A$ be the subgroup generated by the non-involutions of $G$. By (5.1), $A$ is abelian. If $b \in G \setminus A$, then, as $A \cap bA = \emptyset$, $bA$ consists of involutions. Hence for each $a \in A$ we have $ba \cdot ba = 1$, whence $bab = a^{-1}$.

We have $A \subseteq C_G(A) \subseteq G$. Then $\overline{G} = G/C_G(A)$ acts on $A$ faithfully. By the above, each element of $\overline{G}$ must act on $A$ by inversion. It follows that $|G/C_G(A)| = 2$. Since each $b \in C_G(A) \setminus A$ acts on $A$ trivially and simultaneously be inversion, we have $A = C_G(A)$. Thus $|G/A| = 2$ and so $G = A \rtimes Z_2$. □

**Lemma 6.** *Let $G_1$ and $G_2$ be finite groups satisfying (5.1). Suppose that every abelian section of $G_2$ is isomorphic to a section of $G_1$ and that $G_1$ is nonabelian. Then $G_2$ itself is isomorphic to a section of $G_1$.*

Proof. We may assume that $G_2$ is also nonabelian. By Lemma 5, $G_i = A_i \rtimes B_i$ where $A_i$ is abelian and $B_i \simeq Z_2$ acts on $A_i$ by inversion. So it suffices to show that $A_2$ is isomorphic to a section of $A_1$, i.e. that for each prime $p$ the $p$-component of $A_2$ is isomorphic to a section of the $p$-component of $A_1$. For $p$ odd this follows directly from the hypothesis.

Let $T_i$ be the 2-component of $A_i$. Since $T_2$ is abelian, we have

$$T_2 \simeq H/K, \quad \text{where} \quad K \triangleleft H \subseteq G_1 = A_1 \rtimes B_1.$$

Note that every subgroup $X$ of $G_1$ not containing in $A_1$ has the form $X = (X \cap A_1) \rtimes \langle x \rangle$, where $x$ is *any* element from $X \setminus A_1$. Therefore if $K \not\subseteq A_1$, then $T_2 \simeq H/K \simeq (H \cap A_1)/(K \cap A_1)$ is a section of $A_1$ and so of $T_1$. If $K \subseteq A_1$ and $H \subseteq A_1$, then $T_2$ is again a section of $A_1$. Finally, if $K \subseteq A_1$ but $H \not\subseteq A_1$, then $T_2 \simeq ((H \cap A_1)/K) \rtimes \langle h \rangle = \overline{A} \rtimes \langle h \rangle$, where $\langle h \rangle \simeq Z_2$ acts on $\overline{A}$ by inversion. Since $T_2$ is abelian, we have $\exp T_2 = 2$. But then $\langle T_2, B_2 \rangle = T_2 \times B_2$ is abelian and hence is isomorphic to a section of $G_1 = A_1 \rtimes B_1$. It follows that $T_2$ must be isomorphic to a section of $T_1$. □

2. Our aim is to write down a set of disjunctive identities of $D_{2m} = Z_m \rtimes Z_2$ such that the only groups they hold in are sections of $D_{2m}$. The first two d-identities are $\omega_{2m}$ and (5.1). Any group in which they are valid is a finite group of the form $A \rtimes Z_2$, or abelian (Lemma 5). Next, we will add to $\omega_{2m}$ and (5.1) several new d-identities of $D_{2m}$ such that the only *abelian* groups all these d-identities hold in are sections of $D_{2m}$. Then by Lemma 6, *any* group in which all these d-identities hold will be a section of $D_{2m}$.

We assume that *m is even*. The case of odd $m$ is simpler and will be discussed at the end of the section.

**Lemma 7.** *Let $A$ be a finite abelian group of an even exponent $m$. If $A$ has no sections of the form*

$$Z_2^3, \quad Z_4 \times Z_2, \quad Z_p^2, \quad Z_p \times Z_2^2,$$

*where $p$ is an odd prime divisor of $m$, then $A$ is isomorphic to a section of $D_{2m}$.*

P r o o f. Let $A = A_2 \times A_{p_1} \times \ldots \times A_{p_s}$ be the decomposition of $A$ into primary components. Since $A$ has no sections of the form $Z_p^2$, where $p$ is an odd divisor of $m$, each $A_{p_i}$ must be cyclic. Hence $A = A_2 \times Z_k$, where $k$ is an odd divisor of $m$. If $A_2$ is cyclic, then $A$ is a section of $Z_m$, otherwise it follows from the hypothesis that $A \simeq Z_2 \times Z_2$, which is a section of $D_{2m}$. □

3. Since $m$ is even, the identity
$$x^m = 1 \tag{5.2}$$
is valid in $D_{2m}$. Every group satisfying (5.2) has exponent $m$, and so, for each of the four groups of Lemma 7, it remains to find a d-identity of $D_{2m}$ that is not valid in this group.

(i) *The group $Z_p^2$.* Let $m = p^k m_p$, where $p$ is an odd prime and $(p, m_p) = 1$. Consider the formula
$$(x_1^2 = 1) \vee (x_2^2 = 1) \vee (x_1^{m_p} \in \langle x_2^{m_p} \rangle) \vee (x_2^{m_p} \in \langle x_1^{m_p} \rangle), \tag{5.3}$$
where $\alpha \in \langle \beta \rangle$ denotes
$$(\alpha = \beta) \vee (\alpha = \beta^2) \vee \ldots \vee (\alpha = \beta^{p^k - 1}).$$

It is valid in $D_{2m} = Z_m \rtimes Z_2$ because if $x_1^2 \neq 1$ and $x_2^2 \neq 1$, then $x_i \in Z_m = Z_{p^k m_p}$. Hence $x_i^{m_p} \in Z_{p^k}$ and so one of the subgroups $\langle x_i^{m_p} \rangle$ is contained in the other. On the other hand, (5.2) fails in $Z_p^2$: take a basis of this group for the $x_i$.

(ii) *The group $Z_4 \times Z_2$.* Take the d-identity (5.3) with $m_p = m_2$ (that is, $m = 2^k m_2$ with $m_2$ odd). Again it is valid in $D_{2m}$. But it is not valid in $Z_4 \times Z_2$: if $a$ and $b$ are generators for $Z_4$ and $Z_2$, respectively, and $x_1 = (a, b)$, $x_2 = (a, 1)$, then
$$(x_1^2 \neq 1) \,\&\, (x_2^2 \neq 1) \,\&\, (x_1^{m_2} \notin \langle x_2^{m_2} \rangle) \,\&\, (x_2^{m_2} \notin \langle x_2^{m_2} \rangle).$$

(iii) *The group $Z_2^3$.* This group is eliminated by the following generalization of the d-identity (2.4)
$$[\bigvee_{i \neq j} x_i^{m_2} \in \langle x_j^{m_2} \rangle]$$
$$\vee \, [\bigvee_{\sigma \in S_3} x_{\sigma(1)}^{m_2} \in \langle (x_{\sigma(2)} x_{\sigma(3)})^{m_2} \rangle] \vee [\bigvee_{\sigma \in S_3} x_{\sigma(1)}^{m_2} \ni \langle (x_{\sigma(2)} x_{\sigma(3)})^{m_2} \rangle] \tag{5.4}$$
$$\vee \, [\bigvee_{\sigma \in S_3} (x_{\sigma(1)} x_{\sigma(2)})^{m_2} \in \langle (x_{\sigma(2)} x_{\sigma(3)})^{m_2} \rangle]$$
using the same arguments as in §2.

(iv) *The group $Z_p \times Z_2^2$* (where $p$ is an odd divisor of $m$). Let $m = 2^{k_0} p^{k_1} p_2^{k_2} \ldots p_s^{k_s}$, then $m_2 = p^{k_1} p_2^{k_2} \ldots p_s^{k_s}$, $m_p = 2^{k_0} p_2^{k_2} \ldots p_s^{k_s}$. Consider the following formula in three variables $x_1, x_2, x_3$:

$$\left((x_1^{m_p})^{x_2} = x_1^{-m_p}\right) \vee \left((x_1^{m_p})^{x_3} = x_1^{-m_p}\right)$$
$$\vee \left(x_2^{m_2} \in \langle x_3^{m_2} \rangle\right) \vee \left(x_3^{m_2} \in \langle x_2^{m_2} \rangle\right). \tag{5.5}$$

It is valid in $D_{2m} = Z_m \leftthreetimes Z_2$. Indeed, since $m$ is even, we have $\exp D_{2m} = m$ and so $x^{m_p}$ is a $p$-element and $x^{m_2}$ is a 2-element for any $x \in D_{2m}$. If the first two clauses of (5.4) are false, then both $x_2$ and $x_3$ belong to $Z_m$. But then both $x_2^{m_2}$ and $x_3^{m_2}$ are 2-elements, and the rest is obvious. On the other hand, a natural basis of the group $Z_p \times Z_2 \times Z_2$ refutes the d-identity (5.5) in this group.

Combining all of the above, we obtain

**Theorem 4.** *Let $m$ be even. Then the following formulas form a basis of d-identities of $D_{2m}$:*
*$\omega_{2m}$, (5.1), (5.2), (5.4);*
*(5.3) for all prime divisors $p$ of $m$;*
*(5.5) for all odd prime divisors $p$ of $m$.* □

4. The case of an odd $m$ is simpler. First, instead of (5.2), one has to take the d-identity

$$(x^2 = 1) \vee (x^m = 1). \tag{5.2'}$$

Next, instead of Lemma 7, we notice that if $A$ is a finite abelian group satisfying (5.2'), which is not isomorphic to a section of $D_{2m}$, then $A$ has a section isomorphic to either $Z_2^2$ or $Z_p^2$, where $p$ is an odd prime dividing $m$. The latter is eliminated by (5.3), but the former – by the d-identity

$$(x_1^m = 1) \vee (x_2^m = 1) \vee ((x_1 x_2)^m = 1). \tag{5.6}$$

Thus, a basis of d-identities of the group $D_{2m}$ with an odd $m$ consists of the formulas $\omega_{2m}$, (5.1), (5.2'), (5.3) for all primes dividing $m$, and (5.6).

## 6. Identities of regular representations

Using the method of [PK], outlined in the introduction, we can now easily write down bases for the identities of regular representations of the groups considered in the previous sections. Since this theory is applicable to ordinary representations only, the characteristic of the ground field $K$ is always chosen to be coprime with the order of the group in question.

Recall that if
$$\delta = (f_1 = 1) \vee \ldots \vee (f_n = 1)$$

is a disjunctive identity of a group $G$, then $\delta^*$ denotes the corresponding identity of the regular representation $\text{Reg}_K G$, that is,

$$\delta^* = (f_1^{y_1} - 1) \ldots (f_n^{y_n} - 1).$$

It will be convenient to extend this notation to the numbers of d-identities as well. For example, $(2.10)^*$ will denote the element of $KF$ corresponding to the d-identity $(2.10)$, namely the element

$$(x_1^{4y_1} - 1)(x_2^{4y_2} - 1)((x_1 x_2)^{4y_3} - 1)((x_1 x_2^2)^{4y_4} - 1).$$

**Proposition 5.** *Let* $\text{char} K \neq 2$. *Then the following elements of* $KF$ *form a basis of identities of* $\text{Reg}_K D_8$:

$$\prod_{0 \leq i < j \leq 8} \left( (x_i x_j^{-1})^{y_{ij}} = 1 \right), \tag{6.1}$$

$$x^4 - 1, \tag{6.2}$$

$$(x_1^{2y_1} - 1)(x_2^{2y_2} - 1)(x_3^{2y_3} - 1)((x_1 x_2^{-1})^{y_4} - 1)((x_1 x_3^{-1})^{y_5} - 1)((x_2 x_3^{-1})^{y_6} - 1). \tag{6.3}$$

P r o o f. Clearly these elements are $(2.1)^*$, $(2.2)^*$ and $(2.3)^*$, respectively. Since $(2.1)$–$(2.3)$ form a weak basis of d-identities of $D_8$, it suffices to show that the variety of group representations defined by $(6.1)$–$(6.3)$ is locally finite and ordinary. But this variety is obviously contained in $\omega \mathcal{B}_4$, where $\mathcal{B}_4$ is the Burnside variety of groups of exponent 4. Since $\mathcal{B}_4$ is locally finite and $\text{char} K \neq 2$, the result follows from Fact 2 in Introduction. □

It is not harder to find bases for the identities of the regular representations of $Q_8$ and $A_4$. For example, we know that the formula $\omega_{12}$, $(2.7)$ and $(2.8)$ form a weak basis of d-identities of $A_4$. To find a basis of identities of $\text{Reg}_K A_4$ (provided $\text{char} K \neq 2, 3$), we take $\omega_{12}^*$, $(2.7)^*$, $(2.8)^*$ and also the identity $x^6 - 1$. The latter defines the variety $\omega \mathcal{B}_6$, and since $\mathcal{B}_6$ is locally finite, we are done.

**Proposition 6.** *Let* $\text{char} K \neq 2, 3$. *Then a basis of identities of* $\text{Reg}_K S_4$ *consists of* $\omega_{24}^*$, $(2.10)^*$–$(2.13)^*$, *plus the identities*

$$x^{12} - 1, \tag{6.4}$$

$$[[[x_1, x_2], [x_3, x_4]], [[x_5, x_6], [x_7, x_8]]] - 1. \tag{6.5}$$

P r o o f. We know that $\omega_{24}$ and $(2.10)$–$(2.13)$ form a basis of d-identities of $S_4$. Since $S_4$ is a 3-solvable group of exponent 12, the identities $(6.4)$ and $(6.5)$ are valid in $\text{Reg}_K S_4$. Moreover, these two identities define the variety of representations $\omega \mathfrak{V}$, where $\mathfrak{V}$ is the variety of 3-solvable groups of exponent 12. Since $\mathfrak{V}$ is locally finite and $\text{char } K \nmid 12$, the result follows. □

The argument above generalizes to arbitrary solvable groups. Indeed, let $G$ be a finite solvable group such that *char $K \nmid |G|$*, and suppose that $\delta_1, \delta_2, \ldots, \delta_n$ is a (weak) basis of d-identities of $G$. To find a basis of identities of $\text{Reg}_K G$, we take all $\delta_i^*$ and add to them

$$x^e - 1 \quad \text{and} \quad v_s - 1 \tag{6.6}$$

where $e = \exp G$, $s$ is the solvability class of $G$, and $v_s \in F$ is the identity of solvability of class $s$. Then the variety of group representations defined by (6.6) is locally finite and ordinary, and we are done.

For example, if $G = D_{2m}$, then the identities (6.6) are $x^{2m} - 1$ (or $x^m - 1$ if $m$ is even) and $[[x_1, x_2], [x_3, x_4]] - 1$.

If the group in question is not solvable, there are other ways to obtain identities that guarantee that the corresponding variety is locally finite and ordinary. Let

$$s_n(x_1, \ldots, x_n) = \sum_{\sigma \in S_n} (-1)^\sigma x_{\sigma(1)} \ldots x_{\sigma(n)}$$

be the standard polynomial of degree $n$.

**Lemma 8.** *If $G$ is a group of order $n$ and exponent $e$, then $x^e - 1$ and $s_{n+1}$ are identities of $\text{Reg}_K G$. Furthermore, the variety of group representations defined by these two identities is locally finite (and ordinary if $\text{char} K \nmid e$).*

P r o o f. Since $s_{n+1} \equiv 0$ in any algebra of dimension $n$, the first statement is obvious. To prove the second, it is enough to show that if $\rho : G \longrightarrow \text{Aut}_K V$ is a faithful representation satisfying $x^e - 1$ and $s_{n+1}$, then $G$ is locally finite. We may assume that $G \subset \text{End}_K V$, and then $s_{n+1}(g_1, \ldots, g_{n+1}) = 0$ for any $g \in G$. Let $A$ be the subalgebra of $\text{End}_K V$ generated by $G$. Since $s_{n+1}(x_1, \ldots, x_{n+1})$ is a multilinear polynomial, it is clear that $s_{n+1}(x_1, \ldots, x_{n+1}) \equiv 0$ is identically true on $A$. Hence $G$ is a periodic subgroup of a $PI$-ring $A$. It is known that such a group is locally finite (see, for example, [Pr]). $\square$

Thus, if we know a basis (or at least a weak basis) of d-identities of a group $G$, the procedure of finding a basis of identities of $\text{Reg}_K G$ becomes automatic. For example, combining Theorem 3, Lemma 8 and Fact 2, we obtain

**Proposition 7.** *Let $\text{char } K \neq 2, 3, 5$. Then the identities (4.1)*–(4.7)* plus $x^{60} - 1$ and $s_{361}$ form a basis of identities of $\text{Reg}_K A_6$.* $\square$

School of Mathematics, Institute for Advanced Study, Princeton, NJ 08540
*E-mail address*: vovsi@math.ias.edu, vovsi@math.rutgers.edu